\documentclass[12pt]{article}

\usepackage[T1,T2A]{fontenc}
\usepackage[utf8]{inputenc}

\usepackage{amsfonts}
\usepackage{amssymb}
\usepackage{amsmath}
\usepackage{amsthm}

\def \le {\leqslant}
\def \ge {\geqslant}

\usepackage{tikz} 
\usetikzlibrary{calc,arrows,decorations.pathreplacing,fadings,3d,positioning}

\begin{document}

\begin{Large}
\centerline{A note on Dirichlet spectrum}
\end{Large}
\vskip+1cm
\centerline{ by {\bf Renat K. Akhunzhanov} and {\bf
Nikolay G. Moshchevitin}\footnote{Research is supported by the Russian Science Foundation under grant 19-11-00001.
 }
}
\vskip+1cm

{\bf 1. Structure of Dirichlet spectrum.}
\vskip+0.3cm

Let   $ \frak{g}(\pmb{y}), \pmb{y}\in \mathbb{R}^n$ be an arbitrary norm in 
$\mathbb{R}^n$. Dirichlet spectrum 
$\mathbb{D}_{\frak{g}}^{[n]} $ for  simultaneous approximation with respect to the norm $\frak{g}(\cdot)$ is defined as follows.
For $\pmb{\theta} = (\theta_1,...,\theta_n)\in \mathbb{R}^n$ consider the  irrationality measure function
$$
\psi_{\frak{g},\pmb{\theta}} (t) = 
\min_{q\in \mathbb{Z}_+:\, q\le t }\,\,\,\,
\min_{\pmb{p} = (p_1,...,p_n)\in \mathbb{Z}^n}\,\,\,\,
\frak{g} (q\pmb{\theta} - \pmb{p})
$$
and define 
$$
d^{[n]}_{\frak{g}}(\pmb{\theta}) = \limsup_{t\to \infty} t (\psi_{g,\pmb{\theta}} (t))^n.
$$
Then 
$$
\mathbb{D}_{\frak{g}}^{[n]}  = \{ d\in \mathbb{R}:\,\,
\exists \, \pmb{\theta}\in \mathbb{R}^n \setminus\mathbb{Q}^n\,\,
\text{such that}\,\, d =  d_{\frak{g}}^{[n]}(\pmb{\theta}) \}.
$$
We should note that if we
consider the critical determinant  
$
\frak{K}_{\frak{g}}
$
of the cylinder
$$
\{ \pmb{z} = (x, \pmb{y})\in \mathbb{R}^{n+1}:\,\,\,
\frak{g}(\pmb{y})\le 1, |x|\le 1\}
$$
(for the definitions see \cite{Cas}) then
$$
\mathbb{D}_{\frak{g}}^{[n]} \subset
\left[0, \frac{1}{ \frak{K}_{\frak{g}}}
\right].
$$
In particular it follows 
from
Minkowski's convex body theorem   that 
$$
 \frak{K}_{\frak{g}}\ge 
 \frac{\Omega[\frak{g}]}{2^n}
,\,\,\,\,\text{where}\,\,\,\,
\Omega[\frak{g}] =\int_{\pmb{x}\in \mathbb{R}^n: \, \frak{g}(\pmb{x}) \le 1} \, d\pmb{x} .
$$
So
$$
\mathbb{D}_{\frak{g}}^{[n]} \subset
\left[0, \frac{2^n}{\Omega[\frak{g}]}
\right] .
$$

For  $n=1 $ and the standard norm $\frak{g}(y) = |y|$ the spectrum $ \mathbb{D} = \mathbb{D}^{[1]}_{|\cdot |}$
was studied by many authors (\cite{Sz,SDav, DiN, Di,  Iva, Ivaa, Iva1},  see also some recent results concerning 
related metric settings for uniform approximation  to one number \cite{Kl,Kl1,Kl2}).
In the case $n\ge 2$ not much is known, however a complete structure of the spectrum
$\mathbb{D}^{[2]}_{|\cdot |} $ for Euclidean norm
$$
|\pmb{y}| = \sqrt{y_1^2+y_2^2}
$$
in $ \mathbb{R}^2$ was discovered by Akhunzhanov and Shatskov in \cite{AS}. It turned out that 
$
\mathbb{D}^{[2]}_{|\cdot |}
$ 
is a segment, namely
$$
\mathbb{D}^{[2]}_{|\cdot |} = \left[ 0, \frac{2}{\sqrt{3}}
\right].
$$ 
Here the value of $ \max  \mathbb{D}^{[2]}_{|\cdot |} =  \frac{2}{\sqrt{3}}$ is related to the critical determinant of the cylinder
$$
\{ (x, y_1,y_2)\in \mathbb{R}^3:\,\,\, y_1^2+y_2^2\le 1, \, |x|\le 1\}
$$
calculated by Mahler \cite{M}.

As far as the authors know,
$
\mathbb{D}^{[2]}_{|\cdot |}
$ 
is the only 
Diophantine spectrum with completely known structure.
The authors believe that  for $n\ge 2$ and arbitrary norm $\frak{g}(\cdot)$ the spectrum
$\mathbb{D}_{\frak{g}}^{[n]} $ should have the same simple structure, namely it should be the segment of the form $ 
\left[0,  \frac{1}{\frak{K}_{\frak{g}}}\right]$, however up to now they are not able to prove this even in the case $n=3$ for Euclidean norm in $\mathbb{R}^3$.

We would like to mention that very recently a series of results concerning spectrum
$
\mathbb{D}^{[2]}_{\frak{g}}
$
for arbitrary norm $\frak{g}$ in $\mathbb{R}^2$  was obtained in \cite{KA}.
In particular
it was proven there that the maximal point  $
\max \, \mathbb{D}^{[2]}_{\frak{g}}
$
of  the spectrum $
\mathbb{D}^{[2]}_{\frak{g}}
$ is not isolated in $
\mathbb{D}^{[2]}_{\frak{g}}
$.

In this paper we are interested in a more detailed analysis of distribution of the values
$d_{\frak{g}}^{[2]}(\pmb{\theta}) = d_{\frak{g}}^{[2]}(\theta_1,\theta_2) $ in the two-dimensional case $\pmb{\theta} =(\theta_1,\theta_2) \in \mathbb{R}^2$.

\vskip+0.3cm

{\bf 2. Badly approximable points, best approximations and Dirichlet improvability.}
\vskip+0.3cm

Let us use the notation $|\pmb{y}|$ for the Euclidean norm of $\pmb{y}\in \mathbb{R}^n$.
A point $\pmb{\theta} =(\theta_1,...,\theta_n)\in \mathbb{R}^n$ is called {\it badly approximable} if there exists a positive $\gamma$ such that
\begin{equation}\label{bad}
\min_{\pmb{p} \in \mathbb{Z}^n}
|q \pmb{\theta} -\pmb{p}|
> \frac{\gamma}{q^{1/n}}\,\,\,\,\, \forall q \in \mathbb{Z}_+,
\end{equation}
or equivalently for any norm $\frak{g}(\cdot) $ in $\mathbb{R}^n$ one has
$$
\inf_{t\ge 1}\, t (\psi_{\frak{g},\pmb{\theta}} (t))^n >0.
$$
We should note that $\psi_{\frak{g},\pmb{\theta}} (t)$ is a piecewise constant function and in the case 
$ \pmb{\theta}\in \mathbb{R}^n \setminus\mathbb{Q}^n$ one has $\psi_{\frak{g},\pmb{\theta}} (t)>0 $ for every $t$.
So we can define the unique infinite sequence of integers 
\begin{equation}\label{best}
q_1 = 1<q_2<....<q_{\nu-1}<q_\nu<...
\end{equation}
such  that
$$
\psi_{\frak{g},\pmb{\theta}} (t) = \psi_{g,\pmb{\theta}} (q_{\nu-1})\,\,\,\,\,
\text{for}\,\,\,\,\,
q_{\nu-1}\le t<q_\nu.
$$
Moreover for all $\nu$ large enough ($\nu\ge \nu_0[g]$) there exists unique 
$\pmb{p}\in \mathbb{Z}^n$ such that 
$$
\psi_{\frak{g},\pmb{\theta}} (q_{\nu}) = \frak{g}(q\pmb{\theta}-\pmb{p}).
$$
Of course the sequence (\ref{best}) of the best approximations depends on the norm $\frak{g}$.
However as all the norms in $\mathbb{R}^n$ are equivalent, from Theorem 1 from the paper \cite{AM}
we know that $\pmb{\theta}\in \mathbb{R}^n$ is badly approximable if and only if the inequality 
$$
\sup_{\nu\in \mathbb{Z}_+} 
 \frac{q_{\nu}}{q_{\nu-1}} < \infty
 $$
 holds
 for the sequence of the best approximations in any fixed norm $\frak{g}$.
 This result has an obvious quantitative form which we formulate here for an arbitrary norm $\frak{g}$.
 Define
 $$
 m_{\frak{g}}(\pmb{\theta}) =  \limsup_{\nu\to \infty} \frac{q_{\nu}}{q_{\nu-1}}.
 $$
 
 \vskip+0.3cm
 {\bf Proposition 1.} {\it If the inequality (\ref{bad}) holds for a certain positive $\gamma$ for all $q$ large enough
then for the best approximations with respect to  norm  $\frak{g}$ one has
 \begin{equation}\label{limsup}
 m_{\frak{g}}(\pmb{\theta})< M 
 \end{equation}
 with 
 $$ M =  c_1(n;\frak{g}) \gamma^{-n}.
 $$ 
 Conversely if (\ref{limsup}) holds for the best approximations with respect to  norm $\frak{g}$  for some $M$
 and for all $n$ large enough then  for all $q$ large enough we have (\ref{bad}) with 
 $$\gamma = c_2(n;\frak{g}) M^{-1}.$$ 
 Here $c_{j}(n;\frak{g}),  j = 1,2$ are explicit constant depending on dimension $n$ and norm $\frak{g}$.
 }
 \vskip+0.3cm

  \vskip+0.3cm
  
  Another interesting phenomenon is related to {\it singularity} and {\it Dirichlet improvability}.
  A vector $\pmb{\theta}\in \mathbb{R}^n $ is called  {\it singular} if
  $d_{\frak{g}}^{[n]} (\pmb{\theta}) =0$. It is clear that this definition does not depend on the norm $\frak{g}$ in $\mathbb{R}^n$.
  For $n=1$ real number $\theta$ is singular if and only if it is rational.
  From another hand, for sup-norm $|\pmb{y}|_\infty = \max_{1\le j \le n} |y_j|$ in $\mathbb{R}^n$ it is clear that
  $\mathbb{D}^{[n]}_{|\cdot |_\infty}\subset [0,1]$ and for almost all $\pmb{\theta}\in \mathbb{R}^n$
  one has $d_{|\cdot |_\infty}^{[n]} (\pmb{\theta}) = 1 $ (see  classical paper by Davenport and Schmidt \cite{DSch} as well as recent   papers \cite{KR, KR1,KA} with many metric results and the references therein). A vector $\pmb{\theta}\in \mathbb{R}^n$ is called {\it Dirichlet improvable} (with respect to sup-norm) if $d_{|\cdot|_\infty}^{[n]}(\pmb{\theta}) < 1$.
In the case $n=1$ H. Davenport and W. Schmidt  showed  \cite{SDav} that  a number $\theta\in \mathbb{R}\setminus \mathbb{Q}$ is Dirichlet improvable if  and only if it is badly approximable. So for $n=1$ a number $\theta$ is Dirichlet improvable if either it is badly approximable or singular.
  
  \vskip+0.3cm
  The following result was obtained in \cite{many}.
  
  \vskip+0.3cm
   {\bf Proposition 2.} {\it
   For $ n \ge  2$, the set
 of Dirichlet improvable vectors $\pmb{\theta}\in \mathbb{R}^n$ 
has continuum many points which are neither  badly approximable nor singular.}
  \vskip+0.3cm
  
  By the way, in \cite{many} the authors mentioned that 
  the method of the paper \cite{AS} uses the theory of best approximations and can be adapted to construct Dirichlet improvable points in $\mathbb{R}^2$ that are not simultaneously singular or badly approximable. In the next section we formulate our main results and in particular explain this phenomenon.
  
    \vskip+0.3cm
    At the end of this section we would like to  formulate a  quantitative statement which immediately follows from the argument of a familiar paper \cite{JRUS} by Jarn\'{\i}k (see also  discussion in Section 4.1  in \cite{Msin} as well as Section 2.6 from \cite{C} and \cite{Ce}).
    
      \vskip+0.3cm
   {\bf Proposition 3.} {\it Suppose that   the  components $\theta_1, \theta_2$ of  a vector $\pmb{\theta}\in \mathbb{R}^2$
   are linearly independent over $\mathbb{Q}$ together with 1. Suppose that 
 $$
 d_{|\cdot |}^{[2]}(\pmb{\theta}) \le \epsilon.
 $$
 Then
 $$
  m_{|\cdot |}(\pmb{\theta})\ge \frac{1}{36 \epsilon^2}.
 $$
}
   \vskip+0.3cm

 \vskip+0.3cm

{\bf 3. Main results.}
\vskip+0.3cm
Here we formulate our main results dealing with non-badly approximable  and badly approximable cases.

  \vskip+0.3cm
      {\bf Theorem 1.} {\it 
      Let $ \varphi(t)$ be an arbitrary function increasing to $+\infty$ as $ t \to \infty$.
      For any $\lambda\in \mathbb{D}^2_{|\cdot |} = \left[ 0, \frac{2}{\sqrt{3}}
\right] $ there exists $\pmb{\theta} \in \mathbb{R}^2$ such that

 \noindent
      1) $  d^{[2]}_{|\cdot|} (\pmb{\theta}) = \lambda$;
      
          \noindent
      2)  For the sequence (\ref{best}) for $\pmb{\theta}$ one has
      $\frac{q_{\nu+1}}{q_\nu} \ge  \varphi(n)$ for all $n$, and in particular
      $m_{|\cdot |}(\pmb{\theta})=\infty$.
     }
      \vskip+0.3cm
      
      In particular Theorem 1  together with the result from \cite{AM} show that for any $ \lambda\in \mathbb{D}^{[2]}_{|\cdot |}$ there exists
       $\pmb{\theta}\in \mathbb{R}^2$ which is not badly approximable and $d^{[2]}_{|\cdot|} (\pmb{\theta}) = \lambda$.

   \vskip+0.3cm
      {\bf Theorem 2.} {\it  Let $ \varepsilon \in \left(0,\frac{1}{250}\right)$. Then for any $\lambda$ under the condition
      $\varepsilon<\lambda \le \frac{2}{\sqrt{3}}$ there exists $\pmb{\theta}\in \mathbb{R}^2$ such that 
      
      \noindent
      1) $ \lambda-\varepsilon \le d^{[2]}_{|\cdot|} (\pmb{\theta})\le \lambda$;
      
          \noindent
      2) $m_{|\cdot |}(\pmb{\theta}) < {  10^6}\,  \varepsilon^{-2}$.
     }
      \vskip+0.3cm
Of course we do not take care about the optimality of the constant $ 10^6$ in Theorem 2.
However by comparing the result of Theorem 2 with Proposition 3  we see that the order of the upper bound
$O( \varepsilon^{-2})$ from Theorem 2 is optimal.

In particular Theorem 2 shows that $\mathbb{D}^{[2]}_{|\cdot |}= \left[ 0, \frac{2}{\sqrt{3}}
\right] $ is the closure of the set 
$$
\{ 
\lambda:\,\,\, \exists \,\,\text{badly approximable}\,\, \pmb{\theta}\in \mathbb{R}^2\,\,
\text{such that}\,\,
d^{[2]}_{|\cdot|} (\pmb{\theta})=\lambda\}.
$$

 The method of the proofs of Theorems 1,2 rely on the construction from \cite{AS}.
 
 The structure of the paper is as follows.
 In Sections 4 we introduce all necessary parameters. In Section 6 we describe the inductive construction and formulate general Theorem 3. In Section 7 we deduce Theorems 1 and 2  from Theorem 3.
 The rest of the paper (Sections 8 - 11) are devoted to a complete proof of Theorem 3.

 \vskip+0.3cm
{\bf 4. Parameters.}
 \vskip+0.3cm

\vskip+0.3cm
Suppose that $\varepsilon_\nu>0$  form a decreasing sequence
\begin{equation}\label{deka}
\frac{1}{100} \ge \varepsilon_1 \ge  \varepsilon_2\ge ...\ge \varepsilon_\nu\ge \varepsilon_{\nu+1}>...\, .
\end{equation}
We consider a sequence 
of intervals
$$
   \Delta_{\nu} = (\alpha_{\nu},\omega_{\nu}) \subset \left[ 0,\frac{2}{\sqrt{3}}\right] ,\,\,\,\,\, \nu =1,2,3,... 
$$
of lengths $ 4\varepsilon_\nu =\omega_\nu -\alpha_\nu $
and construct the values of $V_\nu $ satisfying certain properties and such that $\frac{V_\nu}{\pi} \in\Delta_{\nu} $.
Instead of intervals $\Delta_\nu$ it is convenient consider subintervals 
$$
   \Delta_{\nu}^* = (\alpha_{\nu}^*,\omega_{\nu}^*) \subset  \Delta_{\nu}. 
$$
of length  ${\varepsilon_\nu}$  such that 
either
            $$
            ( \alpha^*_\nu,\omega^*_\nu) \subset \left[ \varepsilon_\nu,  1-{\varepsilon_\nu}\right] \,\,\,\,\,\,\,\,\,\,
            (\text{{\bf case  1}}^0), 
            $$
or
$$
    ( \alpha^*_\nu,\omega^*_\nu)
\subset \left[1,\frac{2}{\sqrt{3}}-{\varepsilon_\nu}\right]\,\,\,\,\,\,\,\,\,\,
            (\text{{\bf case  2}}^0).
$$
The arguments from the proofs below in {\bf cases 1}$^0$ and {\bf 2}$^0$ differ.
In the sequel we will not write $^*$ for the endpoints $\alpha_\nu^*, \beta_\nu^*$ in order to avoid cumbersome notation.
So in each case we denote the corresponding interval $(\alpha_\nu^*, \beta_\nu^*)$ of length $ {\varepsilon_\nu}$
simply as $(\alpha_\nu, \beta_\nu)$ and refer to the condition of the case. Condition 
\begin{equation}\label{ae}
\varepsilon_\nu\le \alpha_\nu <\omega_n
\end{equation}
gives the inequality
\begin{equation}\label{dva}
\omega_\nu = \alpha_\nu +\varepsilon_\nu\le 2\alpha_\nu.
\end{equation}
Define
\begin{equation}\label{BC}
B_\nu^- =
  \frac{\alpha_{\nu}^2  }{\omega_{\nu-1}\omega_\nu},  \,\,\,\,\,\,\,\,\,    B_\nu^+ =
 \frac{5\omega_{\nu}^2}{\alpha_{\nu-1}\alpha_{\nu}},
\end{equation}
and
  \begin{equation}\label{Hn}
  H_\nu^-= \frac{\alpha_{\nu-1}\alpha_\nu^2}{5\omega_{\nu-1}\omega_\nu^2} ,\,\,\,\,\
H_\nu^+ = \frac{\sqrt{5}\,\omega_{\nu-1}\omega_{\nu}^2}{\alpha_{\nu-1}\alpha_{\nu}^2}.
 \end{equation}
 it follows from (\ref{BC},\ref{Hn}) and  (\ref{ae},\ref{dva}) that
 $$
  B_\nu^-< B_\nu^+,\,\,\,\,\,\,\,
   \frac{1}{40}\le H_\nu^-\le \frac{1}{5}< \sqrt{5} \le H_\nu^+  \le 8\sqrt{5}.
   $$

      We suppose that a sequence of positive integers 
            $
 k_\nu, \,\,\,\, \nu = 1,2,3,...$
for every $n$ satisfies   the condition
\begin{equation}\label{ao}
   K_\nu^-=
\frac{3H_{\nu}^+}{\varepsilon_{\nu}} \le k_\nu \le
K_\nu^+ = \frac{H_{\nu}^-}{4\varepsilon_{\nu}^2}
 \end{equation}
 (the last inequality here follows from       (\ref{ae})).
 Here we should note that 
 \begin{equation}\label{ao1}
 B_\nu^- k_\nu^2 \ge
 \frac{30}{\varepsilon_\nu}\,\,\,\,\,\,\,\,\,\,
 \text{and}
 \,\,\,\,\,\,\,\,\,\,
 \frac{H_\nu^+}{k_\nu} \le \frac{\varepsilon_\nu }{3}.
\end{equation}

 Let us consider for example the situation when
 \begin{equation}\label{con}
 \alpha_\nu = \alpha,\,\,\, \,\,\omega_\nu = \omega, \,\,\, \,\,
 \varepsilon_\nu = \varepsilon_*=\frac{\varepsilon}{4} = \omega-\alpha,\,\,\,\,\, k_n = k
 \end{equation}
 are constant sequences. We explain what restrictions  do we have on our parameters.
We use (\ref{ae}) and (\ref{dva}) to obtain
$$
 B_n^- = \frac{\alpha^2}{\omega^2}    \ge  \frac{1}{4},\,\,\,
 B_n^+ =\frac{5\omega^2}{\alpha^2} \le 20,\,\,\,
 K_n^- \le \frac{24\sqrt{5}}{\varepsilon_*},
 \,\,\,
 K_n^+ \ge \frac{1}{160\varepsilon^{2}_*}.
 $$
 If we take $\varepsilon_* \le10^{-4}$ and choose 
 \begin{equation}\label{cA}
 k =\left[\frac{24\sqrt{5}}{\varepsilon_*}\right]+1
 \end{equation}
 the condition (\ref{ao})  will be satisfied.

 So in the case (\ref{con}) when the sequences of  parameters are constant sequences  and $k$ is chosen as in (\ref{cA})
  we have   the inequality
  \begin{equation}\label{maina}
 B_n^+k_n^2  \le \frac{6\cdot 10^4}{\varepsilon_*^2} <\frac{10^6}{\varepsilon^2}. 
 \end{equation}

\newpage
{\bf 5. Inductive construction.}
 \vskip+0.3cm

For $ Q,R>0$  and $\pmb{v} =(1,v_1,v_2)\in \mathbb{R}^3$ we define the cylinder 
$$
\Pi({\pmb{v}}, Q, R)=
\left\{\pmb{z}=  (x,\pmb{y}) = (x,y_1,y_2)\in\mathbb{R}^{3}:\,\, 0\le x\le  Q,\;
|x{\pmb{v}} - {\pmb{z}}|\leqslant R\right\}.
$$

In this section we describe a variant of a   standard inductive procedure of constructing a sequence
integer points
\begin{equation}\label{po}
\pmb{w}_\nu  =(q_\nu, \pmb{p}_\nu)= (q_\nu,p_{1,\nu}, p_{2,\nu})\in \mathbb{Z}^3,  \,\,\, \pmb{p}_\nu = (p_{1,\nu}, p_{2,\nu})\in \mathbb{Z}^2,\,\,\,\,\,
  \nu=1,2,3,... ,
\end{equation}
corresponding rational points
$$
\pmb{v}_\nu   =(1, \frak{v}_\nu) = \left(1,\frac{p_{1,\nu}}{q_\nu}, \frac{p_{2,\nu}}{q_n}\right)\in \mathbb{Q}^3, \,\,\,\,\frak{v}_\nu = \left(\frac{p_{1,\nu}}{q_\nu}, \frac{p_{2,\nu}}{q_\nu}\right)\in \mathbb{Q}^2,\,\,\,\,\,
  \nu=1,2,3,... ,
$$
  cylinders 
 \begin{equation}\label{cyl}
 \Pi_\nu
 =
 \Pi({\pmb{v}_\nu}, q_\nu, R_\nu) 
 ,
 \,\,\,\,  \,\,\,\,\,\,
  \Pi_\nu^-
 =
 \Pi({\pmb{v}_\nu}, q_{\nu-1} , R_{\nu}^-) 
 ,
  \end{equation}
  where
  $$
  R_\nu =  |q_{\nu-1}\pmb{v}_\nu - \pmb{w}_{\nu-1}| =
  q_{\nu-1}|\frak{v}_\nu - \frak{v}_{\nu-1}| 
  ,\,\,\,\,\,
     R_{\nu}^- = |q_{\nu-2}\pmb{v}_\nu - \pmb{w}_{\nu-2}| = q_{\nu-2}|\frak{v}_\nu - \frak{v}_{\nu-2}| 
  $$
 and extended cylinders 
 \begin{equation}\label{emi}
 \overline{\Pi}_\nu
 =
 \Pi({\pmb{v}_\nu}, q_\nu , R_\nu(1+\varepsilon_\nu)) 
 ,
 \,\,\,\,  
  \overline{\Pi}_\nu^-
 =
 \Pi({\pmb{v}_\nu}, q_{\nu-1} , R_{\nu}^-(1+\varepsilon_{\nu}^-)) 
 ,
 \,\,\,\text {where}\,\,\, \varepsilon_\nu^- = \varepsilon_{\nu-1}^{2}.
 \end{equation}
  It is clear that 
  $$
   {\Pi}_\nu\subset
   \overline{\Pi}_\nu,\,\,\,\,\,\,
    {\Pi}_\nu^-\subset
    \overline{\Pi}_\nu^-.
  $$
\vskip+0.3cm
Our objects  for every $n$ should satisfy the following Conditions 1) - 6).

\vskip+0.3cm

\noindent
Condition 1)  For any $\nu$ vectors $ \pmb{w}_{\nu-2},\pmb{w}_{\nu-1},\pmb{w}_\nu$ form a basis of $\mathbb{Z}^3$.

\vskip+0.3cm
\noindent
Condition 2)   For every $\nu$ we have 
$   B_\nu^- k_\nu^2 \le \frac{q_{\nu}}{q_{\nu-1}} \le B_\nu^+k_\nu^2$ where $B_\nu^{\pm}$ are defined in (\ref{BC}).

\vskip+0.3cm
\noindent
  Condition 3)  
$ 
  \frac{H_\nu^-R_{\nu-1} }{k_\nu}\le
R_\nu   \le 
  \frac{H_\nu^+R_{\nu-1} }{k_\nu},
$ where $H_\nu^{\pm}$ are defined in (\ref{Hn})

 \vskip+0.3cm    
\noindent
Condition 4) $\Pi_\nu  \cap \mathbb{Z}^3 =  \overline{\Pi}_\nu  \cap \mathbb{Z}^3
= 
\{ \pmb{0},\pmb{w}_{\nu-1},\pmb{w}_{\nu },\pmb{w}_{\nu }-\pmb{w}_{\nu-1}\}$.

\vskip+0.3cm 
\noindent
Condition 5) $\Pi_\nu^-  \cap \mathbb{Z}^3 =  \overline{\Pi}_\nu^-  \cap \mathbb{Z}^3
= 
\{ \pmb{0},\pmb{w}_{\nu-2},\pmb{w}_{\nu-1 }\}.
$

\vskip+0.3cm
\noindent
Condition 6)   For the volume  
$
V_\nu  = {\rm vol} \,   \Pi_\nu  =
\pi q_\nu (R_\nu)^2
$
  for every $ \nu =1,2,3,...$ one has
$
V_\nu /\pi \in (\alpha_{\nu},\omega_\nu).
$

 \vskip+0.3cm

   \vskip+0.3cm
      {\bf Theorem 3.} {\it 
      For a given sequence (\ref{deka}) of $\varepsilon_\nu$ and for a sequence of parameters $k_\nu$ satisfying (\ref{ao})  
      there exists a sequence (\ref{po}) of integer points $\pmb{w}_\nu$ such that all the Conditions 1) - 6) are valid.}
         \vskip+0.3cm

 We  will give a proof of   Theorem 3 in Sections 7  -  11. Namely, in Sections 6 - 10 we introduce all necessary objects and constructions and prove all the auxiliary statements and in 
 Section 11 we complete the inductive step.
 In the next section we will  show that from the existence of a sequence of integer points $\pmb{w}_\nu, \nu=1,2,3,...$ satisfying  Conditions 1) - 6) above, Theorems 1 and 2 follows.
 For this purpose we need to  use Conditions 2), 3),  and 5) 6),   only.
 The rest conditions (Conditions 1), and 4)) can be considered as auxiliary conditions for the first ones, however they have clear geometric interpretation. Moreover they can clarify the construction.

 \vskip+0.3cm  
Here we would like to formulate  few more remarks.

First of all it follows from Conditions  2), 3) and inequalities (\ref{ao1}) that 
\begin{equation}\label{besta}
\frac{q_\nu}{q_{\nu-1}}\ge
 \frac{45}{\varepsilon_\nu }\,\,\,\,\,\,\,\,\,\,
\text{and}
\,\,\,\,\,\,\,\,\,\,
\frac{R_\nu}{R_{\nu-1}}\le  \frac{\varepsilon_n}{3}.
\end{equation}
 
 As cylinder $\Pi_\nu$ does not have integer points inside, by  Mahler's theorem on  the critical determinant  \cite{M} we have
 \begin{equation}\label{cree}
 V_\nu \le \frac{2\pi}{\sqrt{3}}.
 \end{equation}

   \vskip+0.3cm
   {\bf 6.  Proof of Theorems 1 and 2.}
   \vskip+0.3cm
   
   Here we deduce Theorems 1 and 2 from Theorem 3.
   
   Because of  (\ref{besta})  points $\pmb{v}_n$ form a fundamental sequence and the limit point
   $$
   \Theta =(1,\pmb{\theta}) = (1,\theta_1,\theta_2) =\lim_{n\to \infty} \pmb{v}_\nu
   $$
 satisfies
   \begin{equation}\label{limi}
   |\pmb{\theta}-\frak{v}_\nu|=
   |\Theta - \pmb{v}_\nu|
   \le \sum_{j=\nu}^\infty 
   |\frak{v}_j-\frak{v}_{j+1}|
   \le
 \sum_{j=\nu}^\infty \frac{R_{j+1}}{q_{j}}\le 
 \frac{ \varepsilon_{\nu}\cdot R_\nu}{3}
  \sum_{j=\nu}^\infty \frac{1}{q_{j}}\le
   \frac{ \varepsilon_{\nu}\cdot R_\nu}{2q_\nu}
 .
   \end{equation}

   Now we define cylinders
   $$
   \Pi_\nu(\pmb{\theta}) = \Pi ( \Theta, q_\nu, R_\nu (\pmb{\theta})),\,\,\,\,
   \text{where}\,\,\,\,
   R_\nu(\pmb{\theta}) = |q_{\nu-1}\pmb{\theta} - \pmb{p}_{\nu-1}| = q_{\nu-1}|\pmb{\theta} - \frak{v}_{\nu-1}|.
   $$
   It is clear by the definition that     $\pmb{w}_{\nu-1}$ belongs to the
    boundary of  the cylinder $ \Pi_\nu(\pmb{\theta})$. Now we show that 
    \begin{equation}\label{empte}
      \Pi_\nu(\pmb{\theta})\cap \mathbb{Z}^3 =  \{ \pmb{0},\pmb{w}_{\nu-1},\pmb{w}_{\nu }\}.
    \end{equation}
    To prove this, it is enough to show that  
        \begin{equation}\label{empte1}
    \Pi_{\nu-1}(\pmb{\theta})\subset\overline{\Pi}_{\nu}^- 
        \end{equation}
        and that 
          \begin{equation}\label{empte2}
          \pmb{w}_\nu \,\,\,
          \text{belongs to the facet}\,\,\,
\{ x= q_{\nu}\}\,\,\,
\text{of cylinder}\,\,\,
    \Pi_{\nu}
        \end{equation}
        holds for every $\nu$.
        Indeed
        (\ref{empte}) follows from (\ref{empte1},\ref{empte2}) by Condition 5).

        To obtain (\ref{empte1}) 
        by triangle inequality it is enough to prove the inequality
          \begin{equation}\label{empte3}
     q_{\nu-1}|\pmb{\theta} - \frak{v}_{\nu} |+R_{\nu-1}(\pmb{\theta}) \le R_{\nu}^- (1+\varepsilon_{\nu}^-)
    \end{equation}
    (here $   q_{\nu-1}|\pmb{\theta}-\frak{v}_\nu|$ is the distance between the centres of the sections $\{x=q_{\nu-1}\}$ of the cylinders under the consideration
    and we want the  section of  the  cylinder $\Pi_{\nu-1}(\pmb{\theta})$ to be inside the section of  the cylinder $\overline{\Pi}_{\nu}^-$).
        We deduce  from (\ref{limi}) and (\ref{besta}) the inequality
        $$|R_{\nu-1}(\pmb{\theta}) - R_\nu^-| \le q_{\nu-2}|\pmb{\theta}-\frak{v}_\nu|\le
        \frac{q_{\nu-2}}{q_{\nu-1}} \frac{q_{\nu-1}}{q_\nu} q_\nu |\pmb{\theta} -\frak{v}_\nu|\le \frac{\varepsilon_{\nu-1}\varepsilon_\nu^2 R_\nu}{4050}\le 
        \frac{\varepsilon_{\nu-1}^2\varepsilon_\nu^2 R_\nu^-}{\cdot 20250}.
        $$
        Here in the next estimate we take into account that $R_\nu \le \frac{\varepsilon_\nu}{3} R_{\nu-1} \le  \frac{2\varepsilon_\nu}{3} R_\nu^-$.
        From the other hand,
   $$
      q_{\nu-1}|\pmb{\theta}-\frak{v}_{\nu}| \le   \frac{\varepsilon_\nu q_\nu}{45} |\pmb{\theta}-\frak{v}_{\nu}| \le \frac{\varepsilon_{\nu}^2 R_{\nu}}{90} \le \frac{\varepsilon_{\nu}^2 R_{\nu}^-}{225}
      .$$
  So
  $$
    q_{\nu-1}|\pmb{\theta} - \frak{v}_{\nu} |+R_{\nu-1}(\pmb{\theta}) \le R_{\nu}^- \left(1+\frac{2\varepsilon_{\nu-1}^2}{225}\right),
  $$
  and this gives (\ref{empte3}) and (\ref{empte1}).
            As for  (\ref{empte2}),
          it immediately follows form the inequality (\ref{limi}) as for every $n$ we have
         $$
          R_{\nu+1}(\pmb{\theta}) =q_{\nu}|\pmb{\theta}-\frak{v}_{\nu}| 
          <\frac{\varepsilon_\nu
          R_{\nu}(\pmb{\theta})}{4}<  R_{\nu}(\pmb{\theta}) .
        $$

          Relation (\ref{empte})  means that    vectors $\pmb{w}_n$ form the sequence of all best simultaneous approximations to $\pmb{\theta}$.
              So to prove Theorems 1,2 it is enough to choose the corresponding values of parameters and to  understand what is
      $$
      \limsup_{n\to \infty} \frac{
   {\rm Vol}\,    \Pi_\nu(\pmb{\theta}) }{ \pi} .
   $$

     For the volumes
      $
   {\rm Vol}\,    \Pi_\nu(\pmb{\theta}) ={ \pi}  q_\nu R_\nu (\pmb{\theta})^2 $
      of the cylinders  $ \Pi_\nu(\pmb{\theta})$ and (\ref{cree})  we have
   \begin{equation}\label{piV}
      | {\rm Vol}\,    \Pi_\nu(\pmb{\theta})  -V_\nu|
      =\pi q_\nu |R_{n}(\pmb{\theta})^2
      -
      R_{\nu}^2|\le  3\pi q_n R_n |R_{n}(\pmb{\theta})
      -
      R_{\nu}|
      \le \frac{ \varepsilon_\nu^2V_\nu}{30}<  \varepsilon_\nu
   ,
   \end{equation}
   as
   $$
   |R_{\nu}(\pmb{\theta})
      -
      R_{\nu}|
\le
q_{\nu-1}|\pmb{\theta}-\frak{v}_{\nu-1}| \le
\frac{\varepsilon_\nu^2R_\nu}{900} 
   .
   $$

   Let us take $ \lambda \in \left[ 0, \frac{2}{\sqrt{3}}
\right]  $.
   If we take the sequences $\alpha_\nu, \omega_\nu$ converging to $\lambda$ and  $k_\nu$   large enough such that $ 
   B_\nu k_\nu^2
    \ge\varphi (\nu) $ by Conditions 2) and  6)    we get a proof of Theorem 1.
   
   To prove   Theorem 2 we need to take   constant parameters 
   $\alpha_\nu = \lambda-3\varepsilon_*, \omega_\nu = \lambda - 2\varepsilon_*$
   or
   $\alpha_\nu = \lambda-2\varepsilon_*, \omega_\nu = \lambda - \varepsilon_*$
   dependind on {\bf case 1}$^0$ or {\bf 2}$^0$
   with $\varepsilon_* =\frac{\varepsilon}{4}$
      in (\ref{con})  and with $ k_\nu = k$    defined in (\ref{cA}). 
   Then by Condition 6) and  (\ref{piV}) we have
    $  {\rm Vol}\,    \Pi_n(\pmb{\theta})   
    \in \left(
    \lambda - \epsilon, \lambda \right).
    $
    We take into account  Condition 2) and (\ref{maina}),
   and this gives Theorem 2   .

\vskip+0.3cm
{\bf 7. Integer bases and natural coordinates.}
\vskip+0.3cm

Let 
$$\frak{G}=( \pmb{g}_1,\pmb{g}_2,\pmb{g}_3),
\,\,\,
\pmb{g}_j = (q_j,a_{1,j},a_{2,j})\in \mathbb{Z}^3,\,\,\, q_j >0
$$ 
be a basis of the integer lattice $\mathbb{Z}^3$. 
We consider 
the vectors $\pmb{g}_1^*$ , $\pmb{g}_2^*$ and $\pmb{g}_3^*$ defined by
$$
\pmb{g}_1^* = \left(1,\frac{a_{1,1}}{q_1},\frac{a_{2,1}}{q_1}\right)
,
\,\,\,
\pmb{g}_2^* = \left(0,\frac{a_{1,2}}{\sqrt{a_{1,2}^2+a_{2,2}^2}},\frac{a_{2,2}}{\sqrt{a_{1,2}^2+a_{2,2}^2}}\right),\,\,\,
|\pmb{g}_2^* |= 1
$$
and
$$
\pmb{g}_3^* = (0,a,b),\,\,\,\,\,
|\pmb{g}_3^*|=1,\,\,\,\,\,
 \pmb{g}_3^*\,\,
 \text{is orthogonal to}\,\,
 \pmb{g}_2^*\,\,\,\,\,
 \text{and}\,\,\,\,\,
 \pmb{g}_3^* = x\pmb{g}_1+y \pmb{g}_2+z \pmb{g}_3\,\,\,
 \text{with}\,\,\, z>0
$$
which form a basis in $\mathbb{R}^3$.
Coordinates  $(x,y,z) $ in $\mathbb{R}^3$ with respect to the basis 
$\pmb{g}_1^*, \pmb{g}_3^*,\pmb{g}_3^*$ we will call
{\it natural coordinates} associated with $\frak{G}$.
Consider
unit vectors
$$
\pmb{e}_1 = (1,0,0),\,\,\, \pmb{e}_2 = (0,1,0),\,\,\, \pmb{e}_3 = (0,0,1)$$
 and the  unique linear mapping $\mathcal{G}$
such that  $\mathcal{G} \pmb{g}_j^*= \pmb{e}_j,  j =1,2,3$.
 It is clear that $\mathcal{G}$ preserves volume in $\mathbb{R}^3$ and moreover in any
 affine subspace  of the form $\{ x={\rm const}\}$ it preserves Euclidean distances between points.
 We see that 
 \begin{equation}\label{qdh}
 \mathcal{G} \pmb{g}_1 = (q_1,0,0),
 \,\,\,\,\,\,
  \mathcal{G} \pmb{g}_2 = (q_2,d,0),\,\, d>0,
  \,\,\,\,\,
   \mathcal{G} \pmb{g}_3 = (q_3,f,h),\,\, h>0
   \end{equation}
   with some $d,f,h$. Here  in the right hand side of all the equalities we have the natural coordinates of the vectors $\pmb{g}_j$ with respect to $\frak{G}$. It is clear from the orthogonality that 
   \begin{equation}\label{unic}
   q_1 dh = 1.
   \end{equation}
 We  associate with $\frak{G}$ an unimodular  lattice
 $$
 \Gamma_{\frak{G}} = \mathcal{G}\mathbb{Z}^3.
 $$
 It is clear that   vectors $
 \mathcal{G} \pmb{g}_1,  \mathcal{G} \pmb{g}_2 ,
   \mathcal{G} \pmb{g}_3 
   $
   form a basis of 
   $
 \Gamma_{\frak{G}}
 $.

    In particular in  natural coordinates  $(x,y,z)$ with respect 
   to the basis
 \begin{equation}\label{bau}
  \pmb{g}_1=\pmb{w}_\nu,\,\,\,\,\,
 \pmb{g}_2=\pmb{w}_{\nu-1},\,\,\,\,\,
 \pmb{g}_3=
 \pmb{w}_{\nu-2} ,
 \end{equation}
 for the images $ \pmb{w}_j' =  \mathcal{G} \pmb{w}_j$
  we have
   \begin{equation}\label{bau1}
   \pmb{w}_\nu' = (q,0,0),\,\, \pmb{w}_{\nu-1}' = (a_0,d,0),\,\,  \pmb{w}_{\nu-2}' = (g,f,h)
   \,\,\,
\text{with}   
   \,\,\,
   q=q_\nu,\,\, a_0 = q_{\nu-1}, g = q_{\nu-2},\,\, d = R_\nu.
   \end{equation}
   We consider the hyperplane  $\pi_1$ defined by
   $$
   \pi_1 =
   \{ (x,y,z):\,\,\, z= h\},
   $$
   so 
   $w_{\nu-2}\in \pi_1$.
   
   In natural coordinates  cylinders $\Pi_\nu'=\mathcal{G}\Pi_\nu$ and $\Pi_\nu^{-'}=\mathcal{G}\Pi_\nu^{-}$ from (\ref{cyl}) 
   can be defined by
   $$
   \Pi_\nu '= \{ (x,y,z):\,\, 0\le x \le q,\,\, y^2+z^2 \le d^2\},\,\,\,\,\,
   \Pi_\nu^{-'} =
   \{ (x,y,z):\,\, 0\le x \le a_0,\,\, y^2+z^2 \le f^2+h^2\}
   $$
 respectively.
   
   In the sequel 
   we need to consider mapping $F: (x_2,y_2) \mapsto (x_1,y_1)$   introduced in \cite{AS} and defined by
\begin{equation}\label{map}
x_1 =\frac{h(x_2^2+q^2)}{y_2q},
\,\,\,\,\,
y_1 =\frac{hx_2}{q},
\end{equation}
which depends on $q$ and $h$ as parameters. Here we would like to explain the meaning of this mapping.
Consider the unique cylinder
$
\Pi= 
 \Pi (\pmb{v},Q,R)
$
such that the point $ w_{n} = (q,0,0)$ belongs to its boundary and the point $(x_1,y_1, h) \in \pi_1$ is the center of the facet $\{x = Q\}$ of $\Pi$.
So
$$ 
\pmb{v} = \left(
1,\frac{y_1}{x_1}, \frac{h}{x_1}
\right),
\,\,\,\,\,
Q= x_1\,\,\,\,\,\text{and}\,\,\,\,\,
R = \frac{q}{x_1} \sqrt{y_1^2+h^2}.
$$
This cylinder   $\Pi$ can be characterised in a rather different way: cylinder $\Pi$ is the unique cylinder of the form
$\Pi (\pmb{v},Q,R)$
 such that
the center of its facet $\{x = Q\}$  belongs to $\pi_1$, the point   $ w_{n} = (q,0,0)$ belongs to its boundary  and
the line
$$
\{ (x,y,z): \,\,\, x=x_2, \, z=0\}
$$
in the coordinate plane $\{ z=0\}$ is tangent to the boundary of $\Pi$ in a certain point 
$(x_2,y_2,0)$ (for the details see \cite{AS}). The values  $(x_1,y_1)$ and $(x_2,y_2)$ just satisfy the relation 
(\ref{map}).

\vskip+0.3cm
{\bf 8. Ellipses and hyperbolas}
\vskip+0.3cm

We consider coordinate plane $ \mathbb{R}^2 (x,y)$ and  points
$$ 
W = (q,0),\,\,\,\,
A= (a,d),\,\,\,\,\, q,d>0; \,\,\,\,\, Z = (x,y), \,\,\, Z_2 = (x_2,y_2).
$$
We consider  equations
\begin{equation}\label{eq}
(xy_2-yx_2)^2 +(qy)^2 = (qy_2)^2,
\end{equation}
\begin{equation}\label{eqa}
(ay_2-dx_2)^2 +(qd)^2 = (qy_2)^2.
\end{equation}

The following lemma can be verified with  clear direct calculations.
\vskip+0.3cm

{\bf Lemma 1.} {\it 
Let us fix the points $
W = (q,0), Z_2 = (x_2,y_2).
$
Let 
$\frak{E} $ be the $0$-symmetric ellipse  which passes through  the points 
$W, Z_2$, that is 
$ W, Z_2 \in  \frak{E}$.
Suppose that the tangent line to $\frak{E}$ is parallel to the 
$0x$ coordinate axis, so it is of the form
$$
\{ (x,y) \in \mathbb{R}^2:\,\, y =y_2 \}.
$$
Then the points $W, Z_2$ define ellipse $\frak{E} =\frak{E}(q;x_2,y_2)$ uniquely and 
\begin{equation}\label{e}
 \frak{E} =
 \{
 (x,y) \in \mathbb{R}^2:\,\,\, 
 (x,y) \,\,\text{satisfy {\rm  (\ref{eq})}}\,\,
 \}
\end{equation}
}

\vskip+0.3cm
For points  $ 
W = (q,0),
A= (a,d)$
we consider   the branch of  $0$-symmetric  hyperbola
$$
\frak{H} =
\frak{H}  (q;a,d)=
 \{
 (x_2,y_2) \in \mathbb{R}^2:\,\,\,  y_2 \ge 0,\,\,\,
 (x_2,y_2) \,\,\text{satisfy {\rm  (\ref{eqa})}}\,\,
 \}
$$
and the family  $\mathcal{E}(q;a,d)$ of  $0$-symmetric ellipses which pass through the points
$ 
W = (q,0)$ and 
$
A= (a,d)$, that is 
\begin{equation}\label{f}
\
\mathcal{E} 
=
\{
\frak{E}: \,\,   A,W \in \frak{E}
\}.
\end{equation}

{\bf Lemma 2.} {\it 
Let us fix the points $
W = (q,0), A= (a,d).
$ For any 
ellipse $ \frak{E} \in \mathcal{E}$   consider the unique point $ Z_2=(x_2,y_2) = Z_2(\frak{E}), y_2\ge0$ where the tangent line to the ellips $\frak{E}$ is parallel to  $0x$ axis.
Then
$$
\{ (x_2,y_2) :\,\,\, \exists\,  \frak{E} \in \mathcal{E}\,\,
\text{ such that}\,\,
(x_2,y_2) =
 Z_2(\frak{E}) \} =  \frak{H}  (q;a,d).
$$}

Proof. We substitute $(a,d)$ into (\ref{eq}) and obtain (\ref{eqa}).$\Box$

\vskip+0.3cm
Consider the set
$$
\frak{B}_2 (a,d) =
\{
(x_2,y_2) \in \mathbb{R}^2:\,\,
y_2\ge 0,
\,\,\,
(ay_2-dx_2)^2 +(qd)^2 \ge (qy_2)^2
\}
$$
The boundary of the set $\frak{B}_2 (a,d) $
is the union of the coordinate line $0x$ and hyperbola  $\frak{H}  (q;a,d)$.

\vskip+0.3cm
{\bf Lemma 3.} {\it 
Let us take three points
$$
W = (q,0),\,\,\, A= (a,d),\,\,\,  A'= (a+q,d)
$$
and consider two hyperbolas 
$$
\frak{H} =
\frak{H}  (q;a,d)
\,\,\,
\text{and}
\,\,\,
\frak{H}' =
\frak{H}  (q;a+q,d)
.
$$
Then

\noindent
1)  the  intersection
 $
 \frak{H} \cap \frak{H} '
 $
 consist just of one point $  Z'=(x',y')$ and
$x' = \frac{2}{\sqrt{3}} d$;

\noindent
2) for  any $\lambda \in [1, \frac{2}{\sqrt{3}})$
and for curvilinear triangle 
$$
\frak{T}
\subset 
\frak{B}_2 (a,d) \cap
\frak{B}_2 (a+q,d) 
$$
  with vertices $A,A',Z'$ and bounded by the curves 
$\frak{H},\frak{H}'$ and $\ell_1$,
the intersection
$$
 \ell_\lambda
\cap
\frak{T},\,\,\,
\text{where}\,\,\,
\ell_\lambda = 
\{ (x_2, y_2)\in \mathbb{R}^2: \,\, y_2 = d\lambda\}
$$
is a segment $[R_\lambda,S_\lambda ]\subset \ell_\lambda $ of the length $\varepsilon_\lambda q$, where
\begin{equation}\label{var}
\varepsilon_\lambda = \lambda - 2\sqrt{\lambda^2-1}  .
\end{equation}
The endpoints of this segment have coordinates
\begin{equation}\label{end}
R_\lambda=
\left(r, d\lambda\right),\,\, r =  a\lambda+q\sqrt{\lambda^2-1}
,\,\,\,
S_\lambda=
\left( 
 s
,
d\lambda\right),\,\,
s=
(a+q)\lambda-q\sqrt{\lambda^2-1}
;
\end{equation}
\noindent
3)  for any $a_0 \in \mathbb{R}$ the intersection 
\begin{equation}\label{inter}
\frak{B}_2 =
\bigcap_{k\in \mathbb{Z} } \frak{B}_2 (a_0+kq,d) 
\end{equation}
lies in the strip
\begin{equation}\label{sig}
  \Sigma = \left\{ (x,y):\,\, 0\le x \le \frac{2}{\sqrt{3}} d\right\}
  \end{equation}
  and includes the strip
  $$
  \Sigma_0=\left\{ (x,y):\,\, 0\le x \le   d\right\},
  $$
  so
  $
    \Sigma_0
  \subset
\frak{B}_2
\subset
  \Sigma.
  $
}

\vskip+0.3cm

We should note that in the case $\lambda \in \left(1,\frac{2}{\sqrt{3}}\right)$ for $\varepsilon_\lambda$ defined in (\ref{var}) we have bounds
\begin{equation}\label{var1}
3\cdot \left(\frac{2}{\sqrt{3}}-\lambda\right)<
\varepsilon_\lambda \le 1.
\end{equation}
\vskip+0.3cm

Proof of Lemma 3.  Statements 1) and  2) follow from the  direct calculation.
Statement 3) follows from Statement 1) as  for every $ k \in \mathbb{Z}$
the  intersection
 $
\frak{H}  (q;a_0+kq,d)\cap \frak{H}  (q;a_0+(k+1)q,d)
 $
 consist just of one point $  Z_k'=(x'_k,y'_k)$ and
$x'_k = \frac{2}{\sqrt{3}} d$.$\Box$

\vskip+0.3cm

Let parameters   $ \alpha ,\omega$ and $d,v$ be fixed.
We consider the segment 
 \begin{equation}\label{uv1}
 [U^{[1]},V^{[1]}] 
 \,\,\,\,\,
 \text{with}
 \,\,\,\,\,
 U^{[1]}= (v,d\alpha)
,\,\,\,\,\,
V^{[1]} = ( v, d\omega)
\end{equation}
Here we should note that the endpoints of the segment $ [U^{[1]},V^{[1]}] $ do not depend on parameter $\lambda$.
In {\bf case 1}$^0$  when $\alpha, \omega \in (0,1)$ this segment 
belongs to the strip $\Sigma_0$ and so to the set  $\frak{B}_2$. 
The situation  in  {\bf case 2}$^0$ when parameters $\alpha, \omega\in \left(1, \frac{2}{\sqrt{3}}\right)$ is a little bit more complicated.
In the next lemma  for any $\lambda$ we define
a segment  $[U^{[2]}(\lambda),V^{[2]}(\lambda)] \subset \frak{B}_2$
which
 will belong to the segment  
 $
 [U^{[1]},V^{[1]}]
$
with endpoints
 depending on $\lambda$ and
under certain conditions (see formula (\ref{lowseg1}) from the next section).

\vskip+0.3cm
{\bf Lemma 4.}  {\it Suppose that $ a\ge 0$. Let $\lambda \in [1,\frac{2}{\sqrt{3}})$ and
 \begin{equation}\label{u} 
u = \frac{ d v\lambda}{s},\,\,\,
v =\frac{r+s}{2}= \left(a+\frac{q}{2}\right)\lambda
,
 \end{equation}
where $r$ and $s$ are defined in (\ref{end}). 
 Consider the points 
 $$
 V^{[2]}(\lambda) = (v,d) \in [R_\lambda,S_\lambda],
 \,\,\,\,\,\,
 U^{[2]}(\lambda) = (v,u).
 $$
  Then
 $$
 [U^{[2]}(\lambda),V^{[2]}(\lambda)]\subset \frak{B}_2 (a,d) \cap
\frak{B}_2 (a+q,d) 
$$
and the length of the segment $ [U^{[2]}(\lambda),V^{[2]}(\lambda)]$ is equal to $d\omega\left(1-\frac{v}{s}\right)$.}

\vskip+0.3cm
Proof. Both sets 
$
\frak{B}_2 (a,d) $ and
$
\frak{B}_2 (a+q,d)
$ are star bodies with respect to the origin $0$.
So the intersection 
$
\frak{B}_2 (a,d) \cap
\frak{B}_2 (a+q,d)$ is also a star body and 
$$
{\rm conv}\, \{ 0,R_\lambda,S_\lambda \} \subset \frak{B}_2 (a,d) \cap
\frak{B}_2 (a+q,d).
$$
The point $ U^{[2]}(\lambda)$ is just the intersection of the segment $[0,S_\lambda]$ and the line
$\{ (x,y): x = v\}$. As $a\ge 0$ we see that $ [U^{[2]}(\lambda),V^{[2]}(\lambda)] \subset  {\rm conv}\, \{ 0,R_\lambda,S_\lambda\} \subset\frak{B}_2 (a,d) \cap
\frak{B}_2 (a+q,d)$ So the length of  the segment $ [U^{[2]}(\lambda),V^{[2]}(\lambda)] $ is  equal to
$d\lambda \left(1-\frac{v}{s}\right)$ and everything is clear.$\Box$
\vskip+0.3cm
Now we specify a little bit the choice of parameters.

Let
\begin{equation}\label{apa}
0\le a_0 <q,\,\,\,\,\, a_k = a_0 + kq,\,\,\, k \in \mathbb{Z},\,\,\, k \ge 2
\end{equation}
and let
$$
 V^{[2]}_k(\lambda) = (v_k, d\lambda),\,\,\,\,\,\,
U^{[2]}_k(\lambda)  = (v_k, u_k)
$$
be a point with  coordinates $v_k = v, u_k = u$  depending on $\lambda$ and defined by (\ref{end},\ref{u})  by taking $ a = a_k$,
so   in particular 
$$
u_k =  
\frac{ d v_k\lambda}{s_k},\,\,\,
v_k = \left(a_0+ \left(k+\frac{1}{2}\right)q\right)\lambda,\,\,\,
 s_k = (a_0+(k+1)q)\lambda- q\sqrt{\lambda^2-1}
 $$
 and
 \begin{equation}\label{parat}
 \frac{v_k}{s_k}\ge
 \frac{a_0 +\left(k+\frac{1}{2}\right)q}{a_0+(k+1)q} \ge 1 - \frac{1}{2k}.
 \end{equation}
 Here we should note that 
for the length of the segment $[U_k^{[2]}(\lambda),V_k^{[2]}(\lambda)]$ from Lemma 4  by (\ref{parat}) we have an upper bound
 \begin{equation}\label{lowseg}
 {\rm length}\, [U_k^{[2]}(\lambda),V_k^{[2]}(\lambda)]
 = d\lambda \left(
 1-\frac{v_k}{s_k}\right)
 \le  \frac{d\lambda}{2k} .
 \end{equation}
 We consider the segment $  [
  U^{[1]},
V^{[1]}
  ]
  $
  where $v=v_k$ depends on $\lambda$ and is chosen afterwards. Then
  \begin{equation}\label{lowseg1}
  [U_k^{[2]}(\lambda),V_k^{[2]}(\lambda)] 
  \subset 
  [
  U^{[1]},
V^{[1]}
  ]
  \,\,\,\,\,\,\text{provided}
  \,\,\,\,\,\,
  \frac{1}{4k}  < \omega - \alpha.
 \end{equation}

\vskip+0.3cm
{\bf 9.  Ellipses $\frak{E}$ and a two-dimensional lattice.}
\vskip+0.3cm

We consider lattice
\begin{equation}\label{laa}
\Lambda_{q,a_0,d} = \{ 
\pmb{g}_1m+\pmb{g}_2n:\,\,\, m,n \in \mathbb{Z}\},
\,\,\,\text{where}\,\,\,\,
\pmb{g}_1 = (q,0),\,\,\, \pmb{g}_2 = (d,a_0).
\end{equation}
By $ \widehat{\frak{E}}$ we denote the ellipse  $\frak{E}$ together with its interior, so 
$\frak{E}$  is the boundary of  $  \widehat{\frak{E}}$.
For our further consideration we need two ellipses
$\frak{E}_1$ and $\frak{E}_2$ which depend on the parameters  $q, a, d\lambda$

We are interested in dependence on parameter $\lambda \in \left(0,\frac{2}{\sqrt{3}}\right)$.
However we should differ two cases.

In {\bf case 1$^0$}  we consider a segment $[\alpha,\omega ] \subset (0,1)$ and we simply  deal with $\lambda \in [ \alpha,\omega]$.
In  this  case we define 
 \begin{equation}\label{ttt1}
t_1 = \frac{1}{\omega} = 1+\delta_1>1,
\,\,\,\,\,
\delta_1 =  \frac{1}{\omega} -1>0.
\end{equation}
We take arbitrary $a\in \mathbb{R}$ and  for $ \lambda \in [\alpha, \omega]$ we  consider ellipse
$$
\frak{E}_1 =\frak{E}(q; a, d\lambda).
$$
The second ellipse is defined as a dilated ellipse of the form
\begin{equation}\label{e21}
\frak{E}_2 = t _1\cdot \frak{E}(q; a, d\lambda)
\end{equation}
where parameter of dilatation  $t_1$ is defined by (\ref{ttt1}).

 {\bf Case 2$^0$} 
 is a little bit more complicated. We take $ [ \alpha,\omega]\subset  \left(1,\frac{2}{\sqrt{3}}\right)$
  and deal with $\lambda \in [ \alpha,\omega]$.

In {\bf case 2$^0$}  we define ellipses $\frak{E}_1$ and $ \frak{E}_2$  for these values of $\lambda$.
The first ellipse is defined by 
\begin{equation}\label{e1}
\frak{E}_1 =\frak{E}(q; v_k,d\lambda).
\end{equation}
The second ellipse is defined as a dilated ellipse of the form
\begin{equation}\label{e22}
\frak{E}_2 = t _2\cdot \frak{E}(q; v_k,d\lambda),
\end{equation}
where parameter of dilatation  $t_2$ does not depend on $\lambda$ and is defined by the formula
\begin{equation}\label{ttt2}
t_2 = 1+\delta_2,\,\,\,\,\, \delta_2 =\frac{2}{\sqrt{3}}-\omega.
\end{equation} 
We should note that the values of the coefficients of dilatations $t_1$ and $t_2$ do not depend on  $\lambda$. They depend 
on the endpoints of segment 
$ [ \alpha,\omega]$.

\vskip+0.3cm
{\bf Lemma 5.}  {\it In  both {\bf cases 1}$^0$ and {\bf2}$^0$
 for any integer $k\ge 2$ and for  both ellipses $  \widehat{\frak{E}}_j, j = 1,2$ defined above one has
$$
  \widehat{\frak{E}}_j \cap  \Lambda_{q,a_0,d} =
 \{ 0 , \pm \pmb{g}_1\}
 ,
 $$
 where
 $\pmb{g}_1$ is defined in (\ref{laa}).}

\vskip+0.3cm
Proof. 
As $  \widehat{\frak{E}}_1\subset \widehat{\frak{E}}_2$ it is enough to 
proof Lemma 5 for ellipse  $ \widehat{\frak{E}}_2$ only. 

In  the  {\bf case 1$^0$}  we see from
(\ref{ttt1},\ref{e21})  that  ellipse $ \widehat{\frak{E}}_2$ lies in the strip
$\{ (x,y) :\, -d<x<d\}$, and lemma follows immediately.

Let us consider the  {\bf case 2$^0$}.
Elipse $\frak{E}_1$ can be defined by the equation (\ref{eq}) in coordinates $(x,y)$  with parameters
$$
x_2= v_k = \left(a_0+\left(k+\frac{1}{2}\right) q\right)\lambda,\,\,\,\, y_2 = d\lambda,
$$
that is
$$
\left(xd\lambda  - y\left(a_0+\left(k+\frac{1}{2}\right) q\right)\lambda\right)^2 + (qy)^2 = (qd\lambda)^2.
$$
Applying linear transformation
$$
\left(
\begin{array}{c}
x
\cr 
y
\end{array}
\right)
\mapsto
\left(
\begin{array}{c}
x -\tau y
\cr 
y 
\end{array}
\right),\,\,\,\, \tau = \frac{ \left(a_0+\left(k+\frac{1}{2}\right) q\right)\lambda}{d \lambda}
$$
we see that it is enough to consider ellipse $\frak{E}_2$ defined by the equation
$$
(xd\lambda)^2 + (qy)^2 = (qd\lambda)^2
$$
and lattice $\Lambda_{q,\frac{q}{2},d}$.
 The only "dangerous" lattice points for the dilated ellipse $ \widehat{\frak{E}}_2=t_2 \cdot \widehat{\frak{E}}_1$
 are the lattice points 
 $\left(\pm \frac{q}{2}, \pm d\right)\in\Lambda_{q,\frac{q}{2},d}.
 $
  These point do not belong to $ \widehat{\frak{E}}_2$  if and only if
$$
 \left( \left(\frac{a_0}{q} +k\right)\lambda -\left(\frac{a_0}{q}+\left( k+\frac{1}{2}\right)\right)\lambda\right)^2 + 1 >(t_2\lambda)^2,
$$
or
 $
 t_2< \frac{1}{4}+\frac{1}{\lambda^2} \le  \frac{1}{4}+\frac{1}{\lambda^2}$
 The choice of $t_2$ by (\ref{ttt2}) satisfies this condition.$\Box$

\vskip+0.3cm
{\bf 10. Mapping $F$.}
\vskip+0.3cm

We consider mapping $F: (x_2,y_2) \mapsto (x_1,y_1)$ defined by
(\ref{map}) in the end of Section 5 
and introduced in \cite{AS}.

\vskip+0.3cm
{\bf Lemma 6.} {\it Let   $0<\alpha<\omega<1 $  and 
\begin{equation}\label{delt}
v^2 >q^2\cdot \left(\frac{\alpha\omega}{\omega-\alpha} \cdot \frac{d}{h}-1\right).
\end{equation}
Consider the segment $[U^{[1]},V^{[1]}]$  defined in (\ref{uv1}). Then the image
$F([U^{[1]}, V^{[1]}])$
is the segment
$[F(U^{[1]}),F(V^{[1]})]$ parallel to $0y$ axis, its endpoints have coordinates
$$
 F(U^{[1]}) =\left(   \frac{h(v^2+q^2)}{dq\alpha},\frac{hv}{q}\right)
 ,\,\,\,\,
 F(V^{[1]}) =\left(  \frac{h(v^2+q^2)}{dq\omega},\frac{hv}{q}\right)
$$
and for its length we have lower bound
\begin{equation}\label{lg}
{\rm length} \, [F(U^{[1]}),F(V^{[1]})] 
> q
\end{equation}
}
\vskip+0.3cm
Proof. We use formulas (\ref{map}) to see that 
$$
{\rm length} \, [F(U^{[1]}),F(V^{[1]})]=\frac{h(v^2+q^2)}{dq\alpha}- \frac{h(v^2+q^2)}{dq\omega}
\ge\frac{h(v^2+q^2)}{dq} \cdot \frac{\omega-\alpha}{\alpha\omega}
 >q.
$$
Everything is proven.$\Box$

\vskip+0.3cm
{\bf Lemma 7.} {\it Let 
$\lambda \in [\alpha ,\omega]\subset \left(1,\frac{2}{\sqrt{3}}\right)$.
Consider segment
$[U^{[2]}(\lambda), V^{[2]}(\lambda)]$  defined in
 Lemma 4.
Then this image $F([U^{[2]}(\lambda), V^{[2]}(\lambda)])$
is the segment $[F(U^{[2]}(\lambda)), F(V^{[2]}(\lambda))]$
parallel to $0y$ axis and its endpoints have coordinates
\begin{equation}\label{mapuv}
 F(U^{[2]}(\lambda)) =\left(   \frac{h(v^2+q^2)}{dq\lambda}\cdot \frac{s}{v},\frac{hv}{q}\right)
 ,\,\,\,\,
 F(V^{[2]}(\lambda)) =\left(  \frac{h(v^2+q^2)}{dq\lambda} ,\frac{hv}{q}\right)
\end{equation}
Moreover, suppose that
\begin{equation}\label{ke}
a>\frac{3\sqrt{3}}{2\varepsilon_\omega}\cdot q
\end{equation}
and
\begin{equation}\label{ke1}
\frac{h}{d}\ge \frac{\sqrt{3}}{2}.
\end{equation}
Then  for  any $\lambda$ under the consideration the length of the segment  $F[U^{[2]}(\lambda), V^{[2]}(\lambda)]$ we have the lower bound 
\begin{equation}\label{lg}
{\rm length} \, [F(U^{[2]}(\lambda)), F(V^{[2]}(\lambda))]
>q.
\end{equation}
}
\vskip+0.3cm
Proof. Equations (\ref{mapuv}) follow from (\ref{u}) and (\ref{map}).   Then    
$$
{\rm length} \, [F(U^{[2]}(\lambda)), F(V^{[2]}(\lambda))]=  \frac{h(v^2+q^2)}{dq\lambda}\cdot \frac{s}{v}
-\frac{h(v^2+q^2)}{dq\lambda}\ge
 \frac{h(v^2+q^2)}{dq\lambda}\cdot \frac{s-r}{s+r}
.
$$
But $s-r =\varepsilon_\lambda q$ and $ s+r = (2a+q)\lambda$. Now we continue the lower bound for the length of 
$[F(U^{[2]}(\lambda)), F(V^{[2]}(\lambda))] $ applying  the
inequalities $ v\ge r\ge a\omega$, (\ref{ke1}) and the condition (\ref{ke}). So we obtain 
$$
{\rm length} \, [F(U^{[2]}(\lambda)), F(V^{[2]}(\lambda))]\ge \frac{3\sqrt{3}a}{2\varepsilon_\lambda} \ge  \frac{3\sqrt{3}a}{2\varepsilon_\omega}\ge q
.
$$
This is 
the desired lower bound (\ref{lg}).$\Box$

\vskip+0.3cm
{\bf 11. Inductive   step.}

   \vskip+0.3cm
   We suppose that 
   the points $\pmb{w}_j, 1\le j \le \nu$ satisfying Conditions 1) - 6)   from Section 5  are constructed.
   We  must explain how to construct  point $\pmb{w}_{\nu+1}$ such that the objects from Section 5 satisfy the required conditions  for $(\nu+1)$-th step.

   First of all we  will take vector $\pmb{w}_{\nu+1}$ of the form
   $$
   \pmb{w}_{\nu+1} = \pmb{w}_{\nu-2} + m\pmb{w}_{\nu-1} + k\pmb{w}_{\nu},\,\,\,\, m,k \in \mathbb{Z},
   $$
   so $\pmb{w}_{\nu+1}$ belongs to the affine hyperplane 
   $$ \pi_1= \pmb{w}_{\nu-2}+\langle \pmb{w}_{\nu-1}, \pmb{w}_{\nu}\rangle_\mathbb{R}
   $$
    and lies in the lattice $\Lambda_1 =  \pmb{w}_{\nu-2} +  \Lambda_{q,a_0,d} $
   which is congruent to the lattice  $\Lambda_{q,a_0,d}$. We see that Condition 1) is satisfied.
   
   Then we will use natural coordinates  $(x,y,z)$ with respect to 
   to the basis (\ref{bau}) 
 described in Section 6.
 As it was mentioned in Section 6 in  these coordinates for the basis vectors $\pmb{w}_{\nu-2}',\pmb{w}_{\nu-1}',\pmb{w}_\nu'$  we have (\ref{bau1}) and
 for the values  $q,h,d$ we have formula (\ref{unic}).
   Let $\Gamma = \Gamma_{\frak{G}}$
   We define the two-dimensional lattice
$$
\Lambda=\mathbb{Z}^3\cap {\rm span}\, ( \pmb{w}_{\nu-1}',\pmb{w}_\nu')
$$
with basis $  \pmb{w}_{\nu-1}',\pmb{w}_\nu'$ and determinant
$$
{\rm det}_2\,  \Lambda  =\frac{1}{h} = qd.
$$
The lattice $\Lambda $ lies in the coordinate plane $\{z =0\}$. We can identify it with the lattice  $\Lambda_{q, a_0, d}$  defined in (\ref{laa}) in  Section 9.
 Now by (\ref{unic}) we get
\begin{equation}\label{l}
\lambda_- = qd^2 = \frac{d}{h} = \frac{{\rm Vol}\,\Pi_n }{\pi} \in (\alpha_\nu,\omega_\nu)
\end{equation}
by inductive assumption.
   In coordinates $(x,y,z)$ plane $\pi_1$ has equation
   $\{ z= h \}$. The lattice $\Lambda_1 = \pi \cap \Gamma \subset \pi_1$ consists of all the points with coordinates
   $$
   x= f+ lq,\,\,\,
   y = g+md,
   \,\,\,
   z=h, \,\,\,\,\, 
      l,m \in \mathbb{Z}^3.
   $$
  Lattice $\Lambda_1$ splits into one-dimensional lattices  parallel to $ \pmb{w}_\nu' = (q,0,0)$, namely 
  $$
  \Lambda_1 = \bigcup_{m\in \mathbb{Z}} \{ \pmb{w}_{\nu-2}'+ m \pmb{w}_{\nu-1}'+l\pmb{w}_\nu', l \in \mathbb{Z}\}.
  $$
 Each of the one-dimensional lattices $$\{ \pmb{w}_{\nu-2}'+ m \pmb{w}_{\nu-1}'+l\pmb{w}_\nu', l \in \mathbb{Z}\}$$ 
 belongs to the line
 $$ \ell_m =\{ \pmb{w}_{\nu-2}'+ m \pmb{w}_{\nu-1}'+\zeta\pmb{w}_\nu', \zeta \in \mathbb{R}\}$$
 Euclidean distance between  neighbouring lines $ \ell_m$ and $\ell_{m+1}$ is equal to $d$. So these lines can be enumerated as
  \begin{equation}\label{lines}
      \ell_t =
      \{ (x,y,z): \,\, y = td + \eta, z=h \},\,\,\,\, t\in \mathbb{Z}
      \end{equation}
      with certain $ \eta \in \mathbb{R}$.

   Our aim is to find the point $ \pmb{w}_{\nu+1}$ as a  pre-image of a non-zero point from 
   $\Lambda_1$. We will define this lattice point by means of its image under 
   the mapping $ \frak{G}$, that is we explain how to choose $ \pmb{w}' \in \Gamma$ such that 
   $\pmb{w}_{\nu+1} = \frak{G}^{-1} \pmb{w}'$.

   \vskip+0.3cm
To complete inductive step we need the following easy  lemma.

\vskip+0.3cm
{\bf Lemma 8.} {\it 
Suppose that $\xi, d>0$, $\eta \in \mathbb{R}$ and $ \omega >\alpha$.
Let 
\begin{equation}\label{extre}
\xi (\omega - \alpha ) >d.
\end{equation}
Then
there exist $\lambda\in (\alpha, \omega)$ and  an  integer $ t$  
such  that
\begin{equation}\label{ko}
\xi \lambda = dt +\eta
\end{equation}

 }.

\vskip+0.3cm
Proof.
 When $\lambda $  changes in an  interval $J$ of length $\delta =\omega-\alpha>0$ the value of
 $\xi \lambda$ changes in the dilated interval  $\xi  J$ of length $>d$.$\Box$
      \vskip+0.3cm

      Now we consider the next interval $ (\alpha, \omega) = (\alpha_{\nu+1},\omega_{\nu+1})$ from the Condition 8) from the $(\nu+1)$-th step of inductive process.
     We deal with 
        $k_{\nu+1}$ which satisfies (\ref{ao})
        and
         $$
       \xi =  
      \left(\frac{a_0}{q} + \left(k_{\nu+1} +\frac{1}{2}\right)\right) .
      $$
       By  (\ref{Hn}) and  (\ref{ao}) we have
      $$
       \xi =  
      \left(\frac{a_0}{q} + \left(k_{\nu+1} +\frac{1}{2}\right)\right) \lambda 
      \ge
       k_{\nu+1}\ge \frac{3H_{\nu+1}}{\varepsilon_{\nu+1}} \ge
       \frac{3\sqrt{5}}{\varepsilon_{\nu+1}} = 
       \frac{3\sqrt{5}}{\omega_{\nu+1}-\alpha_{\nu+1}} 
               $$
      So we can use  Lemma 8 to define $\lambda$ for $k = k_{\nu+1}$
       to have
      (\ref{ko}) with  $\eta$ from (\ref{lines})
      Now we can consider the value 
           \begin{equation}\label{deduce}
      v= \left(a_0 + \left(k_{\nu+1} +\frac{1}{2}\right)q\right) \lambda =  \left(a+\frac{q}{2}\right) \lambda,\,\,\,\,
   k_{\nu+1} q\lambda 
      \le v\le (k_{\nu+1} +2)q\lambda \le 2k_{\nu+1}q\lambda
     \end{equation}
  which  depends on $k_{\nu+1}$ and  on the chosen value of $\lambda$.
    For some $ t\in \mathbb{Z}$ we have $td+\eta =  \frac{hv}{q}$.
      Then we deal with the segment  
      $
         [F(U^{[1]}),F(V^{[1]})] 
         $
         and the segment 
         $
      [F(U^{[2]}(\lambda)),F(V^{[2]}(\lambda))]
      $
       from Section 9. 
      It is clear that
 $$
   [F(U^{[1]}),F(V^{[1]})]\times \{z=h\} \subset \ell_t  \,\,\, (\text{\bf case 1}^0), \,\,\,\,\,\,
      [F(U^{[2]}(\lambda)),F(V^{[2]}(\lambda))]\times \{z=h\} \subset \ell_t \,\,\,(\text{\bf case 2}).
      $$
 Recall that  the one-dimensional line
      $
      \ell_t
   $
      contains one dimensional  affine lattice
      $$
      \left\{ (x,y,z): \,\, y =  \frac{hv}{q}, z=h \right\} \cap \Lambda_1
      =
      \left\{ (x,y,z): \,\, x = x_0 + lq, \, l \in \mathbb{Z},\,\, y = \frac{hv}{q}, z=h \right\} .
      $$
      This means that each segment of the line $\ell_t$  of  the length $q$ contains a point of the  lattice $\Lambda \subset \Gamma$.
   From 
    (\ref{Hn}) and  (\ref{ao})  we see that 
    $$
    \left(\frac{v}{q}\right)^2
    \ge
    k_{\nu+1}^2q^2 \lambda^2 \ge
    \frac{9 (H_{\nu+1}^+)^2\alpha_{\nu+1}^2}{\varepsilon_{\nu+1}^2}\ge \frac{45\alpha_{\nu+1}}{\varepsilon_{\nu+1}} \ge
    \frac{\alpha_{\nu+1}\omega_{\nu+1}\lambda_-}{\omega_{\nu+1}-\alpha_{\nu+1}} -1,
      $$
      and condition ({\ref{delt}) of Lemma 6 is satisfied. By the same argument
      $$
      \frac{a}{q} \ge k_{\nu+1} \ge \frac{1}{\varepsilon_{\nu+1}} \ge
      \frac{1}{\frac{2}{\sqrt{3}}-\omega_{\nu+1}} \ge
      \frac{3\sqrt{3}}{2\varepsilon_{\omega_{\nu+1}}},
      $$
      because of $ \frac{2}{\sqrt{3}} -\omega_{\nu+1}\ge \varepsilon_{\nu+1}$. So condition (\ref{ke}) of Lemma 7 is satisfied.
      As for the condition (\ref{ke1}), it follows from  (\ref{l}).
  So all the  conditions 
      of Lemmas 6 and 7 are satisfied.      By Lemmas 6 and 7    we see that the segments
      $  [F(U^{[1]}),F(V^{[1]})]\times \{z=h\} \subset \ell $  and    $[ [F(U^{[2]}(\lambda)),F(V^{[2]}(\lambda))]\times \{z=h\} \subset \ell $ 
      have  length $ >q$.
      We see that in {\bf case 1}$^0$ segment $  [F(U^{[1]}),F(V^{[1]})]\times \{z=h\} $ 
      and 
            in  {\bf case 2}$^0$ segment $[ [F(U^{[2]}(\lambda)),F(V^{[2]}(\lambda))]\times \{z=h\}$
      has an integer point 
      $ \pmb{w}'\in \Lambda_1 \subset \Gamma$.
     This  integer point 
      is just what we need and it defines the next point $ \pmb{w}_{\nu+1} = \frak{G}^{-1} \pmb{w}'$ of our inductive process.
     
     Here we should note that $y$-coordinate of the point $\pmb{w}'$ is positive and equal to $ \frac{hv}{q}$. This will help us to establish  Condition 6).

   \vskip+0.3cm
   Up to now we have verified only Condition 1) of $(\nu+1)$-th step. Now we verify all other conditions.
   
    \vskip+0.3cm
   To check Condition 2) we take into account that 
   \begin{equation}\label{stri}
 \pmb{w}' =\pmb{w}_{\nu+1}'=  \left(q_{\nu+1},  \frac{hv}{q}, h\right) \in  [F(U^{[j]}),F(V^{[j]})] \times \{z=h\} ,\,\,\,\, \text{where}\,\,\,j=1 \,\,\text{or}\,\,2
   \end{equation}
   with  $v$ depending on  the chosen $\lambda$ and satisfying (\ref{deduce}).
   However, in the {\bf case 2}$^0$  by (\ref{lowseg1}) we know that 
   $$
    [F(U^{[2]}(\lambda)),F(V^{[2]}(\lambda))]\subset
   [F(U^{[1]})),F(V^{[1]})]
   $$
   (the last segment here is defined for  the special value of $v$ depending on $\lambda$).
   So  in all the cases 
     \begin{equation}\label{stri}
   \pmb{w}' =\pmb{w}_{\nu+1}'
   \in 
     [F(U^{[1]}),F(V^{[1]})] \times \{z=h\},
   \end{equation}
  and we evaluate $q_{\nu+1}$ by means of endpoints of the segments
   $[F(U^{[1]}),F(V^{[1]})]$ 
   by  taking into account the inequality from (\ref{deduce}).  
   
   Recall that that $ q= q_\nu$.
    We use the formulas for the endpoints of the segment    $[F(U^{[1]}),F(V^{[1]})]$
   and obtain the bounds
      $$
   B_{\nu+1}^- k_{\nu+1}^2q_\nu = \frac{\alpha_{\nu+1}^2 k_{\nu+1}^2 q_\nu}{\omega_\nu\omega_{\nu+1}}\le \frac{k_{\nu+1}^2 q_\nu^2\lambda^2}{q\lambda_-\omega_{\nu+1}}\le
   \frac{v^2}{q\lambda_-}\le
   \frac{h(v^2+q^2)}{dq\omega_{\nu+1}}\le q_{\nu+1}\le
   $$
   $$
   \le
      \frac{h(v^2+q^2)}{dq\alpha_{n+1}}\le
       \frac{4k_{\nu+1}^2q^2\lambda^2+q^2}{q\alpha_{\nu+1}\lambda_-}\le
       \frac{5k_{\nu+1}^2\omega_{\nu+1}^2q_\nu}{\alpha_\nu\alpha_{\nu+1}}
          =B^+_{\nu+1}k_{\nu+1}^2 q_\nu
 $$
 (in the upper bound here we use that $ k_{\nu+1}\omega_{\nu+1}\ge  1$).

         We have checked the inequalities from Condition 2) in all the cases.  
 
  \vskip+0.3cm
 
 Now we verify  Condition 3)
 Let us consider points $\pmb{w}' = \frak{G}\pmb{w}_{\nu+1} = (q_{n+1}, y_0, z_0), y_0 = \frac{hv}{q}, z_0 = h$ and $  \frak{G} \pmb{w}_\nu$ and the corresponding points $\frak{v}_{\nu+1}' = \left(\frac{y_0}{q_{\nu+1}}, \frac{z_0}{q_{\nu+1}} \right) $ and $\frak{v}_n'  = \left(0, 0 \right)$ 
in $\mathbb{R}^2$. As $ | \frak{v}_{\nu+1}' - \frak{v}_{1}'|  = | \frak{v}_{\nu+1} - \frak{v}_{1}|$, we omit $'$  in the notation of the points $ \frak{v}_\nu' , \frak{v}_{\nu+1}' $.
 It is clear that 
 for 
 $$
  | \frak{v}_{\nu+1} - \frak{v}_{\nu}|=
 \frac{1}{q_{\nu+1} }\sqrt{y_0^2 + z_0^2}=
 \frac{1}{q_{\nu+1} }\sqrt{\left(\frac{hv}{q}\right)^2 + h^2} =  \frac{d}{q_{\nu+1}\lambda_- }\sqrt{\left(\frac{v}{q}\right)^2 + 1}
 =\frac{R_\nu}{q_{\nu+1}\lambda_- }\sqrt{\left(\frac{v}{q}\right)^2 + 1}
 $$
 we have upper and lowed bounds
 $$
   \frac{\alpha_{\nu+1}R_\nu k_{\nu+1}}{\omega_\nu q_{\nu+1}}
   \le
 | \frak{v}_{\nu+1} - \frak{v}_{\nu}| 
  \le
  \frac{\sqrt{5}\omega_{\nu+1}R_\nu k_{\nu+1}}{\alpha_\nu q_{\nu+1}}
 $$
 as $ k_{\nu+1}\omega_{\nu+1}\ge 1$.
 We take into account Condition 2)  to get
 $$
\frac{H_{\nu+1}^- R_\nu}{k_{\nu+1}}
 =
   \frac{\alpha_{\nu+1}R_\nu }{\omega_\nu B^+_{\nu+1} k_{\nu+1}}
   \le R_{\nu+1} = 
 q_\nu
  | \frak{v}_{\nu+1} - \frak{v}_{\nu}|\le
   \frac{\sqrt{5}\omega_{\nu+1}R_\nu  }{\alpha_\nu B^-_{\nu+1} k_{\nu+1}} = 
    \frac{H_{\nu+1}^+R_\nu}{ k_{\nu+1}}.
 $$
 This gives us Condition 3).

          \vskip+0.3cm
      
     Condition 4) is satisfied by the construction. 
     
     We will use natural coordinates and consider cylinders
     $
       {\Pi}_{\nu+1}'=\frak{G} {\Pi}_{\nu+1},
     \overline{\Pi}_{\nu+1}'=\frak{G} \overline{\Pi}_{\nu+1}
     $.
     By Condition 3) according to (\ref{besta}) and by (\ref{l}) we have
     $$
     R_{\nu+1} \le \frac{R_\nu}{5} = \frac{d}{5} \le \frac{2h}{5\sqrt{3}}<h.
     $$
     So
  the intersections  $\Pi_{\nu+1}'\cap \Gamma$ and
     $\overline{\Pi}_{\nu+1}'\cap \Gamma$ 
     are covered by the planes $ \{ z=0\}$ and $\{ z=h\}$.
     By the construction the center of the facet 
     $\{ x= q_{\nu+1}\}$ of cylinders $\Pi_{\nu+1}'$ and $\overline{\Pi}_{\nu+1}' $ is the point $
     \pmb{w}' = (q_{\nu+1} ,y_0, z_0) $ and the point $F^{-1} (q_{\nu+1},y_0) = (v_k, d\lambda)$ is just the point of the section
     $ \Pi_{\nu+1}' \cap  \{ z=0\}$ where the tangent line is parallel to $0x$ axis (see the explanation of the geometrical meaning of mapping $F$ 
     in the end of Section 7). Then 
     \begin{equation}\label{inkl}
     \Pi_{\nu+1}' \cap  \{ z=0\} \subset \widehat{\frak{E}}_1,
     \,\,\,\,\,
     \overline{ \Pi}_{\nu+1}' \cap  \{ z=0\} \subset \widehat{\frak{E}}_2
    \end{equation}
     where ellipses are defined in Section 8. We explain these inclusions. Ellipse  $\frak{E}_1$  belongs to the boundary of $\Pi_{\nu+1}'$ and this explains the first inclusion in (\ref{inkl}). From the definition   $\frak{E}_2 = t_j\frak{E}_1 , j =1,2$,
     where the coefficients of dilatations  $t_1, t_2$  are different in {\bf cases 1}$^0$ and {\bf 2}$^0$.
    These coefficients are defined in Section   9 in  formulas (\ref{ttt1}) and (\ref{ttt2}).
   The coefficient of dilatation  for cylinder $\Pi_{\nu+1}'$  is  defined  in (\ref{emi}) and is equal to  $1+\varepsilon_{\nu+1}$
  In both cases we see that  
     $$
   1+\varepsilon_{\nu+1}^-\le t_j ,\,\,\,\,\, j = 1,2,
     $$
     because
     $$
     t_1 =\frac{1}{\omega_\nu} \ge 1 +{\varepsilon_\nu}\ge
     1+\varepsilon_{\nu+1},\,\,\,\,\,
     t_2 = 1+\frac{2}{\sqrt{3}} - \omega_\nu \ge 1+\varepsilon_{\nu} >1+\varepsilon_{\nu+1},
     $$
    by the conditions defining {\bf case 1}$^0$ and {\bf case 2}$^0$ in the beginning of Section 4.
    So we explained the second inclusion in (\ref{inkl}).
    Then
    we apply Lemma 5 to see that 
     in the plane $ \{ z=0\}$ in the intersection  $\Pi_{\nu+1}'\cap \Gamma$ there are  just two integer points - $\pmb{0}$ and $\pmb{w}_\nu'$. Now we consider ellipses $\widehat{\frak{E}}_j, j = 1,2$ shifted by the vector $\pmb{w}_{\nu+1}'$.
     These shifted  ellipses lie in the plane  $ \{ z=h\}$.
     So 
     in the plane  $ \{ z=h\}$ in the intersection     $ \overline{ \Pi}_{\nu+1}' \cap  \{ z=h\} $ we have two integer points $\pmb{w}_{\nu+1}'$ and
     $\pmb{w}_{\nu+1}'-\pmb{w}_\nu'$ only. This gives Condition 4) for  the cylinders
      $
       {\Pi}_{\nu+1}=\frak{G}^{-1} {\Pi}_{\nu+1}',
     \overline{\Pi}_{\nu+1}=\frak{G}^{-1} \overline{\Pi}_{\nu+1}'
     $.

       \vskip+0.3cm
      
      Next we explain Condition 5).

First of all we will show that
   \begin{equation}\label{inka}
   \overline{\Pi}_{\nu+1}^-\subset \overline{\Pi}_{\nu}.
   \end{equation}
Let us explain how to prove (\ref{inka}). We should note that 
$R_{\nu+1}$  is the distance between the centres of sections $\{ x = q_{\nu}\}$ of the cylinders
 $  \overline{\Pi}_{\nu+1}^-$ and $ \overline{\Pi}_{\nu}$.
   We use  the triangle inequality
    to see that
    $$
     |R_{\nu+1}^--R_{\nu}| \le
 q_{\nu-1}  |\, |\frak{v}_{\nu+1}-\frak{v}_{\nu-1}| -
 |\frak{v}_{\nu}-\frak{v}_{\nu-1}|\,|\le q_{\nu-1} |\frak{v}_{\nu+1}-\frak{v}_{\nu}|.
   $$
  So
   \begin{equation}\label{dife}
    |R_{\nu+1}^--R_{\nu}| \!\le\!
  q_{\nu-1}|\frak{v}_{\nu+1} - \frak{v}_{\nu}|
   \!\le\!
    q_{\nu}|\frak{v}_{\nu+1} - \frak{v}_{\nu}|
 =R_{\nu+1}\le
      \frac{ R_{\nu} \varepsilon_{n}}{3}.
   \end{equation}
 We see that
 \begin{equation}\label{ded0}
 R_{\nu+1}
 +R_{\nu+1}^-(1+\varepsilon_{\nu+1}^-) \le R_{\nu+1}+ R_{\nu}(1+\varepsilon_{\nu+1}^-)  + |R_{\nu+1}^--R_{\nu}| (1+\varepsilon_{\nu+1}^-).
 \end{equation}
  Now substituting (\ref{besta}), (\ref{dife})   into (\ref{ded0}) we get
 $$
  R_{\nu+1}
 +R_{\nu+1}^-(1+\varepsilon_{\nu+1}^-)  \le R_{\nu}\left(1+{\varepsilon_{\nu+1}^-}
 +
 \frac{2\varepsilon_{\nu}}{3}
 (1+{\varepsilon_{\nu+1}^-})
 \right)
< R_{\nu}(1+\varepsilon_{\nu}).
 $$
 So the section $ \{ x=q_{\nu}\}
 $
 of the cylinder $\overline{\Pi}^-_{\nu+1}$
 is inside the section
  $ \{ x=q_{\nu}\}
 $
 of the cylinder $\overline{\Pi}_{\nu}$,
 and this immediately gives (\ref{inka}).
 
 We have Condition 4) for the $(\nu+1)$-th step of the process yet established.
 So from (\ref{inka})  and Condition 4) we see that 
 $$
   \overline{\Pi}_{\nu+1}^-\cap  \mathbb{Z}^3 \subset \{ \pmb{0},\pmb{w}_{\nu-1},\pmb{w}_{\nu},\pmb{w}_{\nu}-\pmb{w}_{\nu-1}\}.
   $$ 
   To complete the proof of Condition 5) for the $(\nu+1)$-th  inductive step we should show that 
   \begin{equation}\label{should}
   \pmb{w}_{\nu}-\pmb{w}_{\nu-1} \not\in    \overline{\Pi}_{\nu+1}^-.
   \end{equation}
   We use natural coordinates and instead of (\ref{should}) we prove
        \begin{equation}\label{shouldp}
   \pmb{w}_{\nu}'-\pmb{w}_{\nu-1}' \not\in    \overline{\Pi}_{\nu+1}^{-'} = \frak{G}   \overline{\Pi}_{\nu+1}^-.
   \end{equation}
   We deduce a lower bound for the distance  $\rho$ between    
  $$ 
   \pmb{w}_{\nu}'-\pmb{w}_{\nu-1}' =(q_\nu-q_{\nu-1},- d, 0)
   $$ 
   and the center
  $$
  \frac{q_\nu-q_{\nu-1}}{q_{\nu+1}}\, \pmb{w}_{\nu+1}' =
   \left( q_\nu-q_{\nu-1}, \frac{hv(q_\nu-q_{\nu-1})}{q_\nu q_{\nu+1}}, \frac{h(q_\nu-q_{\nu-1})}{q_{\nu+1}}\right).
  $$
   of the section
  $ \{ x = q_{\nu}-q_{\nu-1}\}$ of cylinder $\overline{\Pi}_{\nu+1}^{-'} $.
  Here $v$  satisfies (\ref{deduce}).

  By Condition 3) we see that 
  \begin{equation}\label{d1}
  \rho =
  \left|
    \frac{q_\nu-q_{\nu-1}}{q_{\nu+1}}\, \pmb{w}_{\nu+1}' -
     (\pmb{w}_{\nu}'-\pmb{w}_{\nu-1}')
     \right|\ge
    R_\nu\left( 1+\frac{R_{\nu+1}}{2R_\nu}\right)
    \ge
       R_\nu\left( 1+\frac{H_{\nu+1}^-}{2k_{\nu+1}}\right)
       \ge R_\nu\left(1+ 2\varepsilon_{\nu+1}^2\right)
       .
       \end{equation}
       From the other hand,
       $$
           R_{\nu+1}^-\le R_\nu +\frac{q_{\nu-1}}{q_\nu} R_{\nu+1}\le 
           R_\nu\left(1 + \frac{\varepsilon_\nu^2}{135}\right)
           $$
           and so
     \begin{equation}\label{d2}
       R_{\nu+1}^- (1+\varepsilon_{\nu+1}^-)
       \le
         R_\nu\left(1 + \frac{\varepsilon_\nu^2}{135}\right) (1+\varepsilon_\nu^2)\le
         R_\nu\left(1 + 2\varepsilon_\nu^2\right).
            \end{equation}
            Now from (\ref{d1},\ref{d2}) and the upper bound from (\ref{ao}) we deduce 
            $$ 
            R_{\nu+1}^- (1+\varepsilon_\nu^-) \le \rho,
            $$
            and this gives (\ref{shouldp}).

 So  we proved (\ref{shouldp}) and hence (\ref{should}).

      We have already mentioned that $y$-coordinate of the point $\pmb{w}'$ is positive and equal to $ \frac{hv}{q}$.

             \vskip+0.3cm
             
   Let us check  Condition 6).
   By the construction we have  (\ref{stri}). This means that for a certain $\pmb{\xi} = (v,u) \in [U^{[1]}, V^{[1]}]$ we have
   $\pmb{w}'  = F (\pmb{\xi})$. Then for the second coordinate $u$ we have
   \begin{equation}\label{uuu}
   d\alpha_{\nu+1}   \le u  \le d\omega_{\nu+1}.
   \end{equation}
   So for the volume $ {\rm Vol}\, \Pi_{\nu+1}$ we get
   $$
   {\rm Vol}\, \Pi_{\nu+1} =
   {\rm Vol}\, \Pi_{\nu+1}'=  \pi h\,  {\rm Area}\, \frak{E}(q; a, u) = \pi qhu.
   $$
   Here  ellipse $\frak{E}(q; a, u)$  is defined in Lemma 1  in  the beginning of Section 8,
   $ q= q_\nu$ and the value of $a$ is not important for the calculation of the area.
 As the values  $q,h,d$ satisfy  (\ref{unic}), by (\ref{uuu})
  we establish  Condition 6).

 \vskip+0.3cm
 
 \bf Acknowledgements.}
     
       \vskip+0.3cm

       The second named  is a winner of the   “Leader” contest conducted by Theoretical Physics and Mathematics Advancement Foundation “BASIS” and would like to thank the foundation and jury.

\vskip+0.3cm

 \vskip+1cm

Renat Akhunzhanov,

Astrakhan State University,

e-mail:  akhunzha@mail.ru

\vskip+1cm

Nikolay Moshchevitin

Moscow Center of Fundamental and Applied Mathematics

and

Steklov Mathematical Institute

e-mail: moshchevitin@gmail.com

\end{document}